\documentclass[11pt,a4paper]{article}

\usepackage{amssymb}

\renewcommand{\le}{\leqslant} 
\renewcommand{\ge}{\geqslant}

\def\dd{\mathrm{d}}
\def\ee{\mathrm{e}}

\def\E{\mathbb{E}} 
\def\L{\mathbb{L}}

\renewcommand{\P}{\mathbb P}

\def\ca{\mathcal{A}}
\def\cb{\mathcal{A}'}

\def\pa{A}
\def\pb{A'}

\def\fa{\alpha}
\def\fb{\alpha'}
\def\fc{\alpha''}
\def\na{a}
\def\nb{a'}

\def\cx{{\mathcal X}}
\def\ns{\gamma}

\def\ma{A}
\def\mb{A'}

\def\cm{\mathcal{M}}

\newtheorem{theorem}{Theorem}

\newtheorem{lemma}{Lemma} 
\newtheorem{corollary}[lemma]{Corollary}
\newtheorem{definition}[lemma]{Definition}
\newtheorem{proposition}[lemma]{Proposition}
\newtheorem{problem}[lemma]{Problem}

\def\proof#1{\noindent\trivlist\item[\hskip\labelsep{\bf #1}]\ignorespaces}
\def\endproof{\hfill $\square\qquad$\endtrivlist}

\begin{document}

\title{Harmonic moments  of
  non homogeneous\\ branching processes}
\author{Didier Piau}
\date{\em Universit\'e Lyon 1}
\maketitle 

\begin{abstract}
We study the harmonic moments of Galton-Watson processes, possibly non
 homogeneous, with positive values.  Good estimates of these are
 needed to compute unbiased estimators for non canonical branching
 Markov processes, which occur, for instance, in the modeling of the
 polymerase chain reaction. By convexity, the ratio of the harmonic
 mean to the mean is at most $1$. We prove that, for every square
 integrable branching mechanisms, this ratio lies between $1-\pa/k$
 and $1-\pb/k$ for every initial population of size $k>\pa$.  The
 positive constants $\pa$ and $\pb$ such that $\pa\ge\pb$ are explicit
 and depend only on the generation-by-generation branching
 mechanisms. In particular, we do not use the distribution of the
 limit of the classical martingale associated to the Galton-Watson
 process.
Thus, emphasis is put on non asymptotic bounds and on the
 dependence of the harmonic mean upon the size of the initial
 population. In the Bernoulli case, which is relevant for the modeling
 of the polymerase chain reaction, we prove essentially optimal
 bounds that are valid for every initial population $k\ge1$.
Finally, in the general case and for large enough initial populations,
similar techniques yield sharp estimates of the harmonic
 moments of higher degrees.
\end{abstract}

\bigskip

Date: \today.

Abbreviated title:
Branching harmonic moments



\bigskip

\section*{Introduction}

We study the behaviour of the 
harmonic means $1/\E(1/Z_n)$ of Galton-Watson processes
$(Z_n)_{n\ge0}$, 
possibly
non homogeneous, with positive values. A motivation for this
theoretical problem is the construction of unbiased estimators for
samples of branching Markov processes, when the state of an individual
depends on the number of its siblings. An instance, outside the realm
of pure probability, where this
construction is needed, arises  in
the modeling of the polymerase chain reaction by branching processes,
see Sun~(1995).  In this specific case, the offspring of each
individual is $1$ or $2$, the state of the first descendant is
identical to the state of its parent and the state of the other
descendant, if any, is a stochastic function of the state of its
parent. One wishes to estimate, for instance, the mutation rate of the
reaction from a uniform sample of a given generation. Any unbiased
estimator of the state of such a sample requires to compute the
harmonic mean size of the corresponding generation. But there exists
no closed form of these harmonic means, except for small initial
populations and for small numbers of generations.  Since the mean
sizes of the generations of a branching process are well known, the
above problem is usually circumvented by assuming that the initial
population is very large. Then, an averaging effect occurs which
implies, roughly speaking, that the harmonic mean size of a generation
is close to its mean size.
In the context of the polymerase chain reaction, we showed in previous
 papers,
 see Piau~(2004), that this approximation is accurate for surprisingly
 small initial populations, and we provided sharp quantitative
 estimates of the discrepancy between the harmonic mean and the mean, for
 any initial population. These results also proved useful to establish
 rigorous confidence intervals for the estimator of the mutation rate
 of the polymerase chain reaction, see Piau~(2005).  
Our purpose in the present paper is to give the
 exact extent of this approximation phenomenon for general, possibly
 non homogeneous, Galton-Watson processes with positive values. When
 the approximation phenomenon indeed occurs, we quantify it through
 non asymptotic and essentially optimal bounds.

\section{Results}
\label{s.res}
In the following, $(Z_n)_{n\ge0}$ denotes a positive Galton-Watson
process, possibly non homogeneous.  The distribution of this Markov process
with values in $\{1,2,\ldots\}$
is characterized by a sequence $\Xi:=(\xi_n)_{n\ge1}$ of distributions
on $\{1,2,\ldots\}$, as follows. For every $n\ge1$, 
conditionally on the past of
the process, $Z_n$ is the sum of $Z_{n-1}$ random variables of law
$\xi_n$ which are independent of the past.  Assume that each
$\xi_n$ is integrable of mean $\mu_n\ge1$.  Then $Z_n$ is integrable
and, if $\E_k$ denotes the expectation when $Z_0=k$, 
for any positive integer
$k$, 
$$
\E_k(Z_n)=k\,M_n,\quad\mbox{with}\ M_n:=\prod_{i=1}^n\mu_i.
$$
On the other hand, by convexity, 
the sequence of general term $M_n\,\E_k(1/Z_n)$ is
nondecreasing for $n\ge0$. Thus every term is at least
$1/k$. 
Our aim is to provide explicit bounds of the 
harmonic moments, which imply, 
in particular, that $M_n\,\E_k(1/Z_n)$ is close to
$1/k$ when this is so. In other words, we wish to show 
that the sequence of general
term $M_n\,\E_k(1/Z_n)$ is nearly constant.
Indeed, for every fixed $n\ge0$ and
when $k\to\infty$, the law of large numbers implies that
$\E_k(1/Z_n)$ is equivalent to
$1/(k\,\E_1(Z_n))$, whose value is $1/\E_k(Z_n)=1/(k\,M_n)$.
Much more is true, as we show below. To ease the task of the reader, 
we first state the consequence of
our general results, in
the homogeneous case.

\begin{theorem}
\label{t.a}
Assume that $\Xi$ is constant and square integrable. Thus
$\xi_n=\xi$ and $M_n=\mu^n$ where
$\xi$ is square integrable and $\mu_n=:\mu\ge1$ for every $n\ge1$.
Then, there exists a positive constant $\pa$, which depends only on
 $\xi$, such that, for every integer $k>\pa$ and every $n\ge0$,
$$
1/k\le\mu^n\,\E_k(1/Z_n)\le1/(k-\pa).
$$
Assume furthermore that $\mu\neq1$.
There exists a positive constant $\pb$, which depends only on
 $\xi$, such that $\pb\le\pa$ and, for every integer $k>\pb$,
$$
\lim_{n\to\infty}\mu^n\,\E_k(1/Z_n)\ge1/(k-\pb).
$$
\end{theorem}

\subsection{Harmonic moments}
\label{ss.hm}
Theorem~\ref{t.a} is a consequence of a general quantitative result,
stated as theorem~\ref{t.b} below, which deals
with non homogeneous processes. To state and prove this result, we rely on some
specific families of distributions, that we define now.

\begin{definition}
For every $m\ge1$, the generating function $g_m$
of the positive, integer valued, random
variable $L_m$
is such that, for any $t\in[0,1]$,
$$
\E(t^{L_m}):=g_m(t):=t/(m-(m-1)t).
$$
For any positive $c$, the random variable $L_{c,m}$ is such that,  for
any $t\in[0,1]$, 
$$
\E(t^{L_{c,m}}):=g_{c,m}(t):=(g_m(t^c))^{1/c}
=t/(m-(m-1)t^c)^{1/c}.
$$
\end{definition}

Thus, $g_m=g_{1,m}$.
For every $c$, $L_{c,m}\ge1$ almost surely and $\E(L_{c,m})=m$.
When $m=1$, $L_{m}=L_{c,m}=1$ almost surely.
When $m>1$, the distribution of $L_{m}-1$ is geometric and 
the distribution of $(L_{c,m}-1)/c$ is negative binomial.


\begin{definition}
For any positive $c$, let $\ca_c$ denote the set of distributions of 
integrable random variables
$L\ge1$ such that, for any $t\in[0,1]$,
$$
\E(t^L)\le g_{c,m}(t),\quad m:=\E(L).
$$
For any positive $c$, let $\cb_c$ denote the set of distributions of 
integrable random variables
$L\ge1$ such that, for any $t\in[0,1]$,
$$
\E(t^L)\ge g_{c,m}(t),\quad m:=\E(L).
$$
\end{definition}

Note that one compares  the distribution of 
$L$ to distributions of random variables, not a priori integer
valued but with
the same mean.
We are now able to state our main result.

\begin{theorem}
\label{t.b}
{\bf (1)} Let $n\ge1$.
Assume that there exists $c$ such that
$\xi_i\in\ca_c$ for every $i\le n$. 
Then, for every $k>c$,
$$
M_n\,\E_k(1/Z_n)\le1/(k-c).
$$
{\bf (2)}
Assume that $M_n\to\infty$ when $n\to\infty$ and
that there exists $c$ such that
$\xi_i\in\cb_{c}$ for every $i\ge1$. 
Then, for every $k>c$,
$$
\lim_{n\to\infty}M_n\,\E_k(1/Z_n)\ge1/(k-c).
$$
\end{theorem}

Recall that, by convexity, the sequence $M_n\,\E_k(1/Z_n)$ is
nondecreasing, hence the existence of the limit when $n\to\infty$ 
is a general fact.
Assertion (2) becomes false when $M_n$ is allowed
to stay bounded,
or when
one replaces the limit
$n\to\infty$ by a finite $n$ since, for instance, the $n=0$ value is $1/k$.
On the other hand, in practical situations,
the hypothesis that $M_n\to\infty$ is easy to check
since it only involves the first moments of the
generation-by-generation mechanisms.

The restriction to $k>c$ is important as well. As
proposition~\ref{p.kc}  shows, the behaviours of $\E_k(1/Z_n)$ and
$1/M_n$ can be
quite different if $k$ is not large enough.
Proposition~\ref{p.kc} deals with one generation of a branching
process using random variables distributed as $L_{c,m}$, when $m\to\infty$, and
corollary~\ref{ckc} applies this result to the $n$th generation
of a branching
process using random variables distributed as $L_{c,m}$ for a given $m$, 
when $n\to\infty$. 

\begin{definition}
\label{d.d}
Let $Z$ denote a  random variable and $\P^{c,m}_k$ a probability
measure, such that $Z$ is distributed, with respect to
$\P^{c,m}_k$, like the sum of $k$ 
i.i.d.\ copies of the random variable $L_{c,m}$.
\end{definition}

\begin{proposition}
\label{p.kc}
For any $k\le c$, 
$
m\,\E^{c,m}_k(1/Z)\to\infty
$ 
when $m\to\infty$.
\end{proposition}
\begin{corollary}
\label{ckc}
Assume that $c$ is an integer and that
$\xi_n$ is the distribution of $L_{c,m}$, for every $n$.
Hence $M_n=Z_0\,m^n$.
Then the distribution of $Z_n$ coincides with the distribution of the
first generation of the branching process based on $L_{c,m^n}$.
As a consequence, $m^n\,\E^{c,m}_k(1/Z_n)\to\infty$ when $n\to\infty$.
\end{corollary}

Thus $\E_k(1/Z_n)$ when $k\le c$ may decay on a different scale than $1/M_n$, 
see more on this case in section~\ref{s.cri}.
On the other hand, 
theorem~\ref{t.b} describes every square integrable Galton-Watson
process if $k$ is large enough, as the following theorem shows.

\begin{theorem}
\label{t.c}
Any square integrable distribution on $[1,+\infty[$
belongs to $\ca_c$ for $c$ large enough, respectively to
$\cb_c$ for $c$ small enough.
Conversely, any distribution on $[1,+\infty[$ which belongs to $\ca_c$
is square integrable and its variance is at most $c\,m\,(m-1)$, where
$m$ denotes its mean. Likewise, the variance of
any distribution on $[1,+\infty[$ which belongs to $\cb_c$
is, either finite and at least $c\,m\,(m-1)$, or infinite.
\end{theorem}

We shall precise the optimal values of $c$ for some usual distributions.
Finally, theorem~\ref{t.b} indeed describes the behaviour of
$\E_k(1/Z_n)$ when $k$ is large enough, for any square integrable
branching process.

\subsection{Bernoulli case}

We apply theorem \ref{t.b} to the Bernoulli case when the
offspring is always $1$ or $2$.
This case is relevant in the context of
the polymerase chain reaction. Our techniques yield accurate bounds of
$\E_k(1/Z_n)$ for every positive $k$, that is,
even when $k=1$, for instance in the homogeneous case,
see theorem~\ref{t.dd} below. We first state uniform bounds that are
simple consequences of the results of section~\ref{ss.hm}.

\begin{theorem}
\label{t.d}
Let $n\ge0$.
Assume that $\xi_i=(1-x_i)\,\delta_1+x_i\,\delta_2$ with $x_i$ in $[0,1]$
for every $i\le n$. 
Then, for any $k\ge2$,
$$
1/k\le M_n\,\E_k(1/Z_n)\le1/(k-1).
$$
\end{theorem}

In the homogeneous case, one can prove better bounds.
We write $\E^x_k$ for $\E_k$ when 
$\xi_i=(1-x)\,\delta_1+x\,\delta_2$ for every $i\ge1$.
For every $x$ in $(0,1)$, define 
$$
\fc(x):=1-x,
\quad
\fb(x):=(1-x)/(1+x).
$$
Then $0\le\fb\le\fc\le1$ and 
$\fc$ and $\fb$ decrease from  $\fc(0)=\fb(0)=1$
to
$\fc(1)=\fb(1)=0$.

\begin{theorem}
\label{t.dd}
{\bf (1)} 
For any $k\ge1$ and $n\ge0$,
$$
1/k\le(1+x)^n\,\E^{x}_k(1/Z_n)\le1/(k-\fc(x)).
$$
{\bf (2)} For any $k\ge1$,
$$
\lim_{n\to\infty}(1+x)^n\,\E^{x}_k(1/Z_n)\ge1/(k-\fb(x)).
$$
\end{theorem}
These estimates are precise enough to imply the following side result
about the case $k=1$.
\begin{proposition}
There exists no uniform upper bound of
$(1+x)^n\,\E^{x}_1(1/Z_n)$ for
$n\ge0$ and $x\in(0,1)$, since
$
\lim_{x\to0}\lim_{n\to\infty}(1+x)^n\,\E^{x}_1(1/Z_n)
$
is infinite.
More precisely, for every $x$ in $(0,1)$ and $n\ge1$,
$$
c(x)/x<\lim_{n\to\infty}(1+x)^n\,\E^{x}_1(1/Z_n)<1/x,
$$
where $c(x):=1-x(1-x)/(1+3x)$ is such that $8/9\le c(x)<1$.
\end{proposition}
In theorem~\ref{t.dd}, the value of $\fb(x)$ stems from the general
construction of section~\ref{ss.hm}, but the value of $\fc(x)$
does not.
In other words, a direct application of section~\ref{ss.hm} to the
Bernoulli case yields $\fa(x)$ instead of $\fc(x)$, with
$$
\fa(x):=-\log(1+x)/\log(1-x).
$$
For every $x$ in $(0,1)$, $\fb(x)<\fc(x)<\fa(x)$.

Theorem~\ref{t.dd} follows from the more general case below.
\begin{theorem}
\label{t.e}
Let $\xi_i=(1-x_i)\,\delta_1+x_i\,\delta_2$
for every $i$. 

{\bf (1)} 
If $x_i\ge x$ for every  $i\le n$, then, for any $k\ge1$,
$$
1/k\le\E_k(1/Z_n)\,\prod_{i=1}^n(1+x_i)\le1/(k-\fc(x)).
$$
{\bf (2)} 
If $x_i\le x$ for every  $i$ and
$\displaystyle\sum_{i\ge1}x_i$ diverges, then, for any $k\ge1$,
$$
\lim_{n\to\infty}\E_k(1/Z_n)\,\prod_{i=1}^n(1+x_i)\
\ge1/(k-\fb(x)).
$$
\end{theorem}

\subsection{A discontinuity result}
\label{ss.disc}
In the Bernoulli case, the functions 
$\fc(x)$ and $\fb(x)$ have a nonzero limit at 
$x\to0^+$, hence the second part of theorem~\ref{t.dd} above shows
that the limit of the normalized harmonic moments
does not always depend continuously on the parameters of the model.
We show in this section that the phenomenon 
is general.
For the sake of simplicity, we deal with the 
homogeneous case. 

Let $\cm$ denote a given
subset of $(1,+\infty)$ such that $1$ is a  limit
point of $\cm$.
Below, the limits when $\mu\to1$ are implicitly restricted to $\mu\in
\cm$.
For each $\mu\in\cm$, let $\xi^{\mu}$ denote a distribution of mean $\mu$.
If $\xi_i=\xi^{\mu}$ for every $i\ge1$, define a function $h_k$ on $\cm$
by
$$
h_k(\mu):=\lim_{n\to\infty}\mu^n\,\E_k(1/Z_n). 
$$
\begin{proposition}
\label{p.f}
Assume that, for each $\mu\in\cm$, there exists $\na(\mu)$ and $\nb(\mu)$
such that $\xi^{\mu}$ belongs to $\ca_{\na(\mu)}$ and to $\cb_{\nb(\mu)}$.
Then
$$
1/(k-\nb_*)\le\liminf_{\mu\to1}h_k(\mu)\le\limsup_{\mu\to1}h_k(\mu)
\le1/(k-\na_*),
$$
where
$\na_*:=\limsup_{\mu\to1}\na(\mu)$ and $\nb_*:=\liminf_{\mu\to1}\nb(\mu)$.
Thus, if $\nb_*$ is positive, 
the function $h_k$ is not continuous at $\mu=1^+$.
\end{proposition}
\begin{theorem}
\label{t.f}
In the homogeneous case,
assume that each $\xi^{\mu}$ is the law of $1+X$, where
the law of $X$ is either binomial or Poisson or geometric.
Then, for every $k\ge1$, 
$h_k$ is discontinuous at $\mu=1$, since
$h_k(1)=1/k$ and
$$
\lim_{\mu\to1}h_k(\mu)=1/(k-1).
$$
If the law of $X$ is geometric, then $h_k(\mu)=1/(k-1)$ for every
$\mu>1$ and $h_k(1)=1/k$.

Likewise, assume that each $\xi^{\mu}$ is the law of $L_{c,\mu}$ for a given
positive integer $c$. Then, $h_k(1)=1/k$ for every $k$. Furthermore,
for every $\mu>1$, $h_k(\mu)=1/(k-c)$ 
if $k>c$ and $h_k(\mu)=+\infty$ if $k\le c$.
\end{theorem}

\subsection{Higher harmonic moments}
\label{ss.hh}
We now state an extension of theorem \ref{t.b} to higher harmonic moments.
Theorem~\ref{t.g} is but a special case of proposition~\ref{p.hi} in
section~\ref{s.h}.

\begin{theorem}
\label{t.g}
{\bf (1)} Let $n\ge1$.
Assume that there exists $c$ such that
$\xi_i\in\ca_c$ for every $i\le n$. 
Then, for every positive integer $r$ and every integer $k>rc$,
$$
M_n^r\,\E_k(1/Z_n^r)\le1/[(k-c)(k-2c)\cdots(k-rc)].
$$
{\bf (2)}
Assume that $M_n\to\infty$ when $n\to\infty$ and
that there exists $c$ such that
$\xi_i\in\cb_{c}$ for every $i$. 
Then, for every positive integer $r$ and every integer $k>rc$,
$$
\lim_{n\to\infty}M_n^r\,\E_k(1/Z_n^r)\ge1/[(k-c)(k-2c)\cdots(k-rc)].
$$
\end{theorem}

\begin{corollary}
\label{c.h}
Let $n\ge1$.
Assume that there exists $c$ such that
$\xi_i\in\ca_c$ for every $i\le n$ and write $\sigma^2_k(1/Z_n)$ for
the variance of $1/Z_n$ when $Z_0=k$. 
Then, for every integer $k>2c$,
$$
M_n^2\,\sigma^2_k(1/Z_n)\le(3c)/[k(k-c)(k-2c)].
$$
If, additionally, there exists $c'$ such that
$\xi_i\in\cb_{c'}$ for every $i\le n$, then
the sequence $$
k^3\,M_n^2\,\sigma^2_k(1/Z_n)
$$
is bounded above and below by finite positive constants, independently
of $n$ and $k$, for large enough values of $k$.
\end{corollary}

\subsection{Related studies}
\label{ss.rs}
As mentioned above, Piau (2004) uses preliminary versions of our
results, 
especially in the Bernoulli
case, which is relevant for the study of the polymerase chain reaction. In this
specific case, we are now able to deal directly with every initial
population $k$, even $k=1$.

Ney and Vidyashankar (2003) give asymptotics of the harmonic moments
of every 
integrable homogeneous branching process starting from $k=1$ particle.
When furthermore $L\log L$ is integrable, their
results specialize as follows, see also Bingham (1988) 
for some classical
facts that are recalled below. 

Let $p_1:=\P(L=1)$, $\mu:=\E(L)$, and 
let $\ns$ denote the Karlin--McGregor exponent of the
distribution of $L$ ($\ns$ is also called the
Schr\"oder constant), defined by the equality
$$
p_1\,\mu^{\ns}=1.
$$
Let $W$ denote the almost sure limit of the nonnegative martingale
$Z_n/\mu^n$. 
The Poincar\'e function is
the Laplace transform $P(s):=\E_1(\exp(-sW))$ 
of the distribution of $W$ when $k=1$, and
solves
Poincar\'e's functional equation
$$
P(\mu s)=f(P(s)).
$$
Three cases may arise.
First, when $r>\ns$, $\E_1(1/Z_n^r)/p_1^n$ converges to a finite positive
limit, whose expression is an integral which 
involves the Schr\"oder function $S$, defined for any $t$ in $[0,1)$, by
$$
S(t):=\lim_{n\to\infty}\E_1(t^{Z_n})/p_1^n.
$$ 
Up to a multiplicative constant, $S$ is the unique finite solution on
$[0,1)$ of Schr\"oder's functional equation
$$
S(f(t))=p_1\,S(t).
$$
Second, when $r=\ns$, $\E_1(1/Z_n^{\ns})/(n\,p_1^n)$ 
converges to a finite positive
limit, whose expression involves Poincar\'e and Schr\"oder functions.
Third, when $r<\ns$, $\mu^{nr}\,\E_1(1/Z_n^r)$ converges to a finite positive
limit. Ney and Vidyashankar provide an expression of the limit in terms of
an integral of Poincar\'e function. One can readily 
check that this limit is  in
fact
$\E_1(1/W^r)$ and that the limit is also an upper bound.

When  $L\log L$ is not integrable, 
the results are similar but one must replace
the normalizations $n\,p_1^n=n/\mu^{\ns}$ when
$r=\ns$ and $1/\mu^{nr}$ when $r<\ns$, by similar expressions which
involve the Seneta-Heyde constants.

Coming back definitely to the $L\log L$ case, we recall 
that the distribution of $W$ has a density $w$ on
the nonnegative real line, such that $w(x)/x^{\ns-1}$ is bounded between
two finite positive constants, when $x\to0$, see Dubuc (1971).

The comparison of our results with those 
recalled above is based on two elementary 
lemmas.

\begin{lemma} 
\label{l.8}
For any distribution $\xi$ in $\ca_c$, $\ns(\xi)\ge1/c$.
For any distribution $\xi$ in $\cb_c$, $\ns(\xi)\le1/c$.
\end{lemma}
In other words (see definition~\ref{d.dd} in section~\ref{s.3}),
$$
\mb(\xi)\le1/\ns(\xi)\le\ma(\xi).
$$
\begin{lemma}
\label{l.9}
For any branching process, $k\ge1$, $n\ge0$, and positive $r$,
$$
k^r\,\E_k(1/Z_n^r)\le\E_1(1/Z_n^{r/k})^k.
$$
\end{lemma}
Corollary~\ref{c.10} is not stated as such in the papers that we
mentioned above but it follows from results that we recalled.
\begin{corollary}
\label{c.10}
For any homogeneous branching process of Schr\"oder exponent $\ns$
and any $k>r/\ns$,
the sequence $\mu^{nr}\E_k(1/Z_n^r)$ is bounded as $n$ varies,
by the finite constant $\E_k(1/W^r)$.
\end{corollary}
An interesting feature of
corollary~\ref{c.10} is that it deals with the entire regime
where such a control  of $\mu^{nr}\E_k(1/Z_n^r)$ may hold, namely,
with every initial population $k>r/\ns(\xi)$.
In other words, when $k\le r/\ns(\xi)$, $\mu^{nr}\E_k(1/Z_n^r)$ is not bounded.
Our upper bounds are restricted to higher values of $k$, namely, to
the regime
$k>r\,\ma(\xi)$.

One could think of recovering the dependence with respect to $k$ from
the results of Ney and
Vidyashankar even when $k\ge2$, starting  from the inequality
$$
\E_k(1/W^r)\le\E_1(1/W^{r/k})^k/k^r.
$$
However, the bounds one gets cannot be optimal for $k\ge2$, since
$$
\E_k^c(1/W)=1/(k-c)<\E_1^c(1/W^{1/k})^k/k.
$$
Furthermore, as stressed by Bingham (1988), the law of $W$, hence the
value of $\E(1/W^r)$, may be explicitly computed only in very specific
cases.  In contrast with every other paper we are aware of, the bounds
we provide are explicit.  The assumptions involve only
elementary, step-by-step, mechanisms of the branching process, that
is, the distributions of the number of descendants at each generation.
Also, we allow for inhomogeneous processes, as long as the reproducing
laws belong uniformly to a given space $\ca_c$, respectively $\cb_c$,
and we make explicit the dependence of the bounds on the initial
population.

The introduction of the family of distributions described by 
$g_{c,m}$ for
integer values of $c$ is hardly new, see Harris (1948) for instance.
A key point is that we use them for
noninteger values of $c$ and as a reference scale of any
square integrable distribution. For instance, the $k=1$ case of
Bernoulli distributions requires to make use of values of $c$ in $(0,1)$.
Although
these distributions do not correspond to a branching process 
for noninteger values of $c$,
they still satisfy a semigroup property, and this property
is sufficient for our purposes.
Finally, our methods do not determine the behaviour of the harmonic moments
of
homogeneous processes whose reproducing law is not square integrable.

\subsection{Plan}

The remainder of the paper is organised as follows.
In section~\ref{s.1}, we reduce the case of branching processes in $\ca_c$
and $\cb_c$ to the case of
well-chosen  distributions
$g_{c,m}$, which we solve in section~\ref{s.2}.
In section~\ref{s.3}, we show that the $\ca_c$ and $\cb_c$ cases
imply the result for every square integrable branching process.
In section~\ref{s.h}, we deal with 
harmonic moments of higher degrees. 
In section~\ref{s.g}, we study thoroughly the Bernoulli case, that is,
the case when the
offspring is $1$ or $2$, sharpening our previous results on this
subject. We provide an algorithm to compute the asymptotic harmonic
moments, up to any accuracy, and we present some simulations and
conjectures about this specific case. 
Section~\ref{s.sbo} is a remark about size-biased offsprings.
Finally, in section~\ref{s.cri}, we briefly explain how to deal with
cases when the asymptotic behaviours
of the harmonic mean and the mean do not coincide.

\section{From $g_{c,m}$ to $\ca_c$ and $\cb_c$}
\label{s.1}

We show that every branching process whose 
branching mechanism uses only
laws in $\ca_c$ can be reduced to the case of $L_{c,m}$ 
for a suitable
$m$, and we solve this case. Similar results hold, as regards
the comparison with 
$\cb_c$.

\subsection{Results}

Lemma~\ref{l.3} describes the semi-group structure of each family
$(g_{c,m})_m$. This is the starting point of our computations.
Corollary~\ref{c.6} is a special case of
corollary~\ref{c.5} and corollary~\ref{c.5}
is a consequence of
lemma~\ref{l.3}.
Corollary~\ref{c.5} uses definition~\ref{d.d} in the introduction.

\begin{lemma}
\label{l.3}
For any positive $c$ and any 
$m\ge1$ and $m'\ge1$, $g_{c,m}\circ g_{c,m'}=g_{c,m''}$ with $m'':=mm'$.
\end{lemma}

\begin{corollary}
\label{c.5}
Let $\varphi$ denote a nonnegative completely monotone function.
For every branching process in $\ca_c$, every $k\ge1$ and $n\ge0$,
$$
\E_k(\varphi(Z_n))\le\E_k^{c,m}(\varphi(Z)),
\quad
\mbox{where}\quad m:=M_n.
$$
For every branching process in $\cb_c$, every $k\ge1$ and $n\ge0$,
$$
\E_k(\varphi(Z_n))\ge\E_k^{c,m}(\varphi(Z)),
\quad
\mbox{where}\quad m:=M_n.
$$
\end{corollary}

Recall that $\varphi$ is completely monotone if and only if its derivatives
are such that $(-1)^i\,\varphi^{(i)}$ is nonnegative for every positive
integer $i$.  Nonnegative completely monotone functions are Laplace
transforms of nonnegative measures on $[0,+\infty)$, see chapter~IV of
Widder~(1948). 

\begin{corollary}
\label{c.6}
For every branching process in $\ca_c$ and every positive
$r$, 
$$
\E_k(1/Z_n^r)\le\E_k^{c,m}(1/Z^r),
\quad
\mbox{where}\quad m:=M_n.
$$
For every branching process in $\cb_c$ and 
every positive
$r$, 
$$
\E_k(1/Z_n^r)\ge\E_k^{c,m}(1/Z^r),
\quad
\mbox{where}\quad m:=M_n.
$$
\end{corollary}

\subsection{Proofs}

\proof{Proof of lemma~\ref{l.3}}
Since each $g_{c,m}$ is the conjugate of $g_m$ by the bijection
$t\mapsto t^c$, the case $c=1$ implies the general case.
When $c=1$, $1/(1-g_m(t))$ is an affine 
function of $1/(1-t)$. By
composition,
$g_m\circ g_{m'}$ is also an affine function of $1/(1-t)$ and it
remains to compute its coefficients to prove the semigroup property.
\endproof

\proof{Proof of corollary~\ref{c.5}}
The representation of completely monotone functions which we recalled
after the statement of the corollary shows that
$$
\varphi(z)=\int_0^1t^z\,\dd\pi(t),
$$
for a given measure $\pi$ on $[0,1]$. Thus $\E_k(\varphi(Z_n))$ is a positive
linear functional of the generating function $\E_k(t^{Z_n})$ of $Z_n$.
The function $\E_k(t^{Z_n})$ 
is the $k$th power of the composition from $i=1$ to $n$ 
of the generating function
of $\xi_i$. When the branching process belongs to $\ca_c$, the
generating function of $\xi_i$ is bounded above by $g_{c,\mu_i}$, thus
the composition is bounded above by the composition of the functions
$g_{c,\mu_i}$, which equals $g_{c,m}$. Finally,
$$
\E_k(\varphi(Z_n))\le\int_0^1g_{c,m}(t)^k\,\dd\pi(t)=\E^{c,m}_k(\varphi(Z)).
$$
The proof of the result for branching processes in $\cb_c$ is similar.
\endproof

\proof{Proof of corollary~\ref{c.6}} 
For every positive $r$,
 $\varphi(z):=1/z^r$ is completely mo\-notone. To see this, choose
 $\dd\pi(t)=(\log 1/t)^{r-1}\,\dd t/(\Gamma(r)\,t)$ in the
 representation of $\varphi$ which we used to prove
 corollary~\ref{c.5}.  \endproof

\section{The case $g_{c,m}$}
\label{s.2}

Our task in this section is to evaluate
the moments of $1/Z$ under the measure $\P_k^{c,m}$. The
cases $k>c$ and $k\le c$ yield different asymptotic behaviours of the
first moment of $1/Z$. We begin with the direct way to deal with
$\P_k^{c,m}$ when $k$ is large enough, 
namely, the computation of factorial moments of $Z$
instead of the usual moments, see proposition~\ref{p.e}.
Starting with  lemma~\ref{l.rep}, which gives a representation formula 
valid for every $k$, we study in depth the first harmonic moment, both in
the small $k$ and large $k$ regimes.
Corollary~\ref{c.6} in section~\ref{s.1}
and lemma~\ref{l.psi} below then 
imply theorem~\ref{t.b}.
Lemma~\ref{l.ppi} deals with the case $k=c$.
Lemma~\ref{l.pti}
provides an alternative formulation of the integral of
lemma~\ref{l.rep}, a formulation that lemma~\ref{l.pi} uses to settle
the case $k<c$.
Proposition~\ref{p.kc}, which concludes the case $k\le c$,
is then an easy consequence.


\subsection{Results}
\label{s.21}

We begin with exact formulas.
Theorem \ref{t.g} in section \ref{ss.hh} is a consequence 
of corollary
\ref{c.e} below.
\begin{proposition}
\label{p.e}
{\bf (1)}
For every nonnegative integer $r$,
$$
\E^{c,m}_k(Z(Z+c)\cdots(Z+rc))=m^{r+1}\,k(k+c)\cdots(k+rc).
$$
{\bf (2)}
For every  real number $r$ such that $k>rc$,
$$
m^r\,\E^{c,m}_k(\Gamma((Z/c)-r)/\Gamma(Z/c))
=
\Gamma((k/c)-r)/\Gamma(k/c).
$$
{\bf (3)}
For instance, for every nonnegative integer $r$ such
that $k>rc$,
$$
m^r\,\E^{c,m}_k(1/[(Z-c)\cdots(Z-rc)])=1/[(k-c)\cdots(k-rc)].
$$
\end{proposition}
\begin{corollary}
\label{c.e}
For every nonnegative integer $r$ such
that $k>rc$,
$$
1/k^r\le 
m^{r}\,\E^{c,m}_k(1/Z^{r})
\le
1/[(k-c)\cdots(k-rc)].
$$
For instance, 
for every $k>c$,
$$ 
1/k\le m\,\E_k^{c,m}(1/Z)<m\,\E_k^{c,m}(1/(Z-c))=1/(k-c).
$$
\end{corollary}
Here is a slight generalization of the $r=1$ 
assertion in corollary~\ref{c.e}.
\begin{proposition}
\label{p.sl}
For every $c>0$, $m\ge1$, $u\ge0$ and positive integer $k>u$,
$$
1/k\le m\,\E_k^{c,m}(1/(Z-u))\le1/(k-\sup\{c,u\}).
$$
\end{proposition}


Proposition~\ref{p.e} and corollary~\ref{c.e} are the results that we
use to settle the case $k>c$
in the rest of the paper. We turn to the evaluation of the exact
harmonic moment of $Z$ with respect to $\P_k^{c,m}$. The results below
are mostly used to deal with the case $k\le c$.

\begin{lemma}
\label{l.rep}
For every positive integer $k$, positive $c$ and $m>1$,
$$
\E_k^{c,m}(1/Z)=G(k/c,m)/c,
\quad
G(u,m):=\int_0^1t^{u-1}\,\dd t/(1+(m-1)t).
$$
Alternatively,
$$
G(u,m)=B_{u,1-u}(1-1/m)/(m-1)^u,
$$
where $B_{u,v}$ denotes the incomplete Beta function of parameters $u$
and $v$, that is, for every $x\in[0,1)$,
$$
B_{u,v}(x):=\int_0^xt^{u-1}(1-t)^{v-1}\,\dd t.
$$
\end{lemma}

\begin{lemma}
\label{l.psi} 
Assume that $u>1$. Then,  
$$
m\,G(u,m)\le
1/(u-1).
$$
The order of this upper bound is exact when $m$ is large, since
the function $(m-1)\,G(u,m)$ increases when $m$ increases and converges to
$1/(u-1)$ when $m\to\infty$.
\end{lemma}

\begin{lemma}
\label{l.ppi}
$G(1,m)=(\log m)/(m-1)$.
\end{lemma}

\begin{lemma}
\label{l.pti}
For any $u$,
$(m-1)^u\,G(u,m)$ is an increasing function of $m\ge1$.
\end{lemma}

\begin{lemma}
\label{l.pi}
Assume that $u<1$. When $m\to\infty$, $(m-1)^u\,G(u,m)$ converges
to 
$c_u:=\pi/\sin(\pi u)$. 
Thus, for any $m>1$,
$$
(m-1)^u\,G(u,m)\le c_u.
$$
On the other hand, for any $m\ge2$, $(m-1)^u\,G(u,m)\ge1/(2u)$. 
\\
Bounds of $c_u$ are, for every $u<1$,  
$c_u\ge\pi$ and $1/(2u(1-u))\le c_u\le1/(u(1-u))$.
\end{lemma}
\begin{corollary}
\label{c.p}
Let $k<c$ and $\ell(k,c):=c/(k(c-k))$. For every $m$,
$$
(m-1)^{k/c}\,\E_k^{c,m}(1/Z)\le\ell(k,c).
$$
The order of this upper bound is exact, since
$$
\lim_{m\to\infty}(m-1)^{k/c}\,\E_k^{c,m}(1/Z)\ge\ell(k,c)/2.
$$
\end{corollary}

\subsection{Proofs}

Lemmas \ref{l.ppi}, \ref{l.pti} and \ref{l.pi} stem from the
definitions.
\proof{Proof of proposition~\ref{p.e}}
{\bf (1)}
For any $|x|<1$ and any positive $y$,
$$
1/(1-x)^y=\sum_{r\ge0}x^r\,\Gamma(y+r)/[\Gamma(y)\,\Gamma(r+1)].
$$
Setting $y=Z/c$ and integrating yields
$$
\E_k^{c,m}(1/(1-x)^{Z/c})=\sum_{r\ge0}
\E_k^{c,m}[\Gamma(r+(Z/c))/\Gamma(Z/c)]\,x^r/\Gamma(r+1).
$$
On the other hand,
$$
\E_k^{c,m}(1/(1-x)^{Z/c})=g_{c,m}(1/(1-x)^{1/c})^k=1/(1-mx)^{k/c}.
$$
Using the expansion of $1/(1-mx)^{k/c}$ given above and equating the
coefficients of the two series yield the result for any
nonnegative integer $r$.

{\bf (2)}
For any positive $y$ and $r$ with $y>r$,
$$
\Gamma(r)\,\Gamma(y-r)/\Gamma(y)=\int_0^1t^{y-r-1}\,(1-t)^{r-1}\,\dd t.
$$
Setting $y=Z/c$ and performing the integration yields
$$
\Gamma(r)\,\E_k^{c,m}[\Gamma((Z/c)-r)/\Gamma(Z/c)]
=
\int_0^1g_{c,m}(t^{1/c})^k\,(1-t)^{r-1}\,\dd t/t^{r+1}.
$$
The change of variables $s:=g_{c,m}(t^{1/c})^c=g_m(t)$ yields
$$
\Gamma(r)\,\E_k^{c,m}[\Gamma((Z/c)-r)/\Gamma(Z/c)]
=
\int_0^1s^{k/c-r-1}\,(1-s)^{r-1}\,\dd s/m^{r},
$$
that is, the desired formula.
\endproof
\proof{Proof of lemma~\ref{l.rep}}
Write $\E_k^{c,m}(1/Z)$ as the integral of $g_{c,m}(t)^k/t$ on
$(0,1)$.
Use the change of variable $t':=g_{c,m}(t)^c$.
This yields the first expression of $G$ in the lemma.
To get the expression of $G$ in terms of incomplete Beta function, 
use the change of variable 
$t':=(m-1)t/(1+(m-1)t)$ in
the first expression of $G$.
\endproof
\proof{Proof of lemma~\ref{l.psi}}
In the first expression of $G$ in lemma~\ref{l.rep}, 
use the fact that $1+(m-1)t$ lies
between $m\,t$ and $m$. Thus, $G(u,m)$ lies between the integral of
$t^{u-2}/m$ and the integral of $t^{u-1}/m$, that is, between
$1/((u-1)m)$ and $1/(um)$.
\endproof



\proof{Proof of corollary~\ref{c.p}}
Use the bound of lemma~\ref{l.pi} by
$c_u/c$ for $u:=k/c$, then the bound of $c_u$ by $1/(u(1-u))$.
This yields the bound for every finite value of $m$.
The limit when $m\to\infty$ is $c_u/c\ge1/(2uc(1-u))=\ell(k,c)/2$.
\endproof

\section{From $\ca_c$ and $\cb_c$ to the general case}
\label{s.3}

In this section, we show that every square integrable branching
process belongs to the set $\ca_c$, respectively to the set $\cb_{c}$, 
for a suitable value
of $c$, we prove theorem~\ref{t.c} and we describe the best
possible constants $c$ of theorem~\ref{t.b} in some specific examples.

\subsection{Comparisons}
Our next proposition is related to theorem~\ref{t.c} and justifies
definition~\ref{d.dd} below.
\begin{proposition}
\label{p.17}
If $c_1\le c_2$ and $m\ge1$, $g_{c_1,m}\le g_{c_2,m}$.
If $c_1<c_2$ and $m>1$,  the distribution of $L_{c_2,m}$ belongs to 
$\ca_{c_2}$ but not to $\ca_{c_1}$ and the distribution of $L_{c_1,m}$
belongs
to $\cb_{c_1}$ but not to $\cb_{c_2}$.
Thus, $(\ca_c)_c$ is a strictly increasing sequence and $(\cb_c)_c$ 
is a strictly decreasing sequence.
\end{proposition}
\begin{definition}
\label{d.dd}
For any square integrable distribution $\xi$ on $[1,+\infty[$, let
$$
\ma(\xi):=\inf\{c>0\,;\,\xi\in\ca_c\},
\quad
\mb(\xi):=\sup\{c>0\,;\,\xi\in\cb_c\}.
$$
\end{definition}

\subsection{Examples}

We now study some specific transformations and examples.
Proposition~\ref{p.tr} follows from the definitions.

\begin{proposition}
\label{p.tr}
For every $\xi$, $\mb(\xi)\le \ma(\xi)$.
The inequality is strict except in two cases:
either $\ma(\xi)=\mb(\xi)=0$, and in that case, 
$\xi$ is a Dirac measure
at $m\ge1$;
or $\ma(\xi)=\mb(\xi)=c$ is positive, and in that case,
$\xi$ is the
distribution of a random variable $L_{c,m}$.
\end{proposition}

\begin{proposition}
\label{p.ex}
{\bf (1)}
If the laws of the independent $1+X$ and $1+X'$ belong to $\cb_c$, the
law of $1+X+X'$ belongs to $\cb_{c}$ as well. This statement
is false when one replaces $\cb_c$ by $\ca_c$. 

{\bf (2)}
If the law of  $1+X$ belongs to $\ca_c$ and if $b$ is positive, the
law of $1+bX$ belongs to $\ca_{cb}$. A similar statement
holds if one replaces $\ca_c$ and $\ca_{cb}$ by $\cb_c$ and $\cb_{cb}$. 

{\bf (3)}
If the law of  $L$ belongs to $\ca_c$ and if 
$L'$ dominates stochastically $L$,
then the law of $L'$ belongs to $\ca_c$ as well.
For instance, if $b$ is nonnegative, 
the law of $L+b$
belongs to $\ca_c$. A similar statement
holds if one replaces $\ca_c$ by $\cb_c$.
\end{proposition}

If $\xi$ is a Dirac measure, $\mb(\xi)=\ma(\xi)=0$.
Other usual cases are as follows.

\begin{proposition}
\label{p.ey}
{\bf (1)}
If $\xi$ is uniform on $\{1,\ldots,n\}$, $\ma(\xi)<1$. More
precisely, $2\,n^{\ma(\xi)}=n+1.$

{\bf (2)}
If $\xi$ is uniform on $\{1,n\}$, $\mb(\xi)=(n-1)/(n+1)$ 
and $2^{\ma(\xi)+1}=n+1$.

{\bf (3)}
If $\xi$ is the law of $1+X$ where $X$ 
is binomial $(n,x)$, $\ma(\xi)<1$. More precisely,
$$(1+x\,n)\,(1-x)^{n\ma(\xi)}=1.$$ 

{\bf (4)}
If $\xi$ is the law of $1+X$ where $X$ is Poisson of mean
$x$, $\ma(\xi)<1$. More
precisely,
$$
\ee^{x\ma(\xi)}=1+x.
$$
\end{proposition}

With the notations of section~\ref{ss.rs},
cases (1) to (4) of proposition \ref{p.ey} are such that
$\ma(\xi)=1/\ns(\xi)>\mb(\xi)$.

To check that the three values $\ma(\xi)$, 
$1/\ns(\xi)$ and $\mb(\xi)$ can indeed be
different,
assume that $\xi:=(1-p)\delta_1+(\delta_2+\delta_3)\,p/2$ with
$p\in(0,1)$. 

Then $\ns(\xi):=-\log(1-p)/\log(1+3p)$.
For $p=1/{5}$, one can check that
the function $t\mapsto\E(t^L)/g_{1/\ns,m}(t)$ has
positive derivatives at $t=0$ and $t=1$. 
Thus some values of this function are greater than  $1$ 
and some are
smaller than $1$. This implies that $\ma(\xi)>1/\ns(\xi)>\mb(\xi)$.

\subsection{Bernoulli case}

\begin{definition}
If $\xi=(1-x)\,\delta_1+x\,\delta_2$,
write $\fa(x)$ for $\ma(\xi)$ and
$\fb(x)$ for $\mb(\xi)$.
\end{definition}

\begin{proposition}
\label{p.be}
For any $x\in(0,1)$, $\fb(x)<\fa(x)<1$, since 
$$
\fb(x)=(1-x)/(1+x),
\quad
(1-x)^{\fa(x)}(1+x)=1.$$ 
Thus, $\fa$ and $\fb$ decrease on $(0,1]$ 
from $\fa(0^+)=\fb(0^+)=1$ to
$\fa(1)=\fb(1)=0$. Both are
discontinuous at $0$ since 
$\fa(0)=\fb(0)=0$.
\end{proposition}

One can note that
$\fb(x)\le1-x\le \fa(x)$.

\subsection{Proofs}

\proof{Proof of proposition~\ref{p.17}}
Compare the logarithmic derivatives.
\endproof
\proof{Proof of theorem~\ref{t.c}}
Both results stem from the expansion of $g_{c,m}$ near $1$, which
reads as follows, when $t=o(1)$,
$$
g_{c,m}(1-t)=1-m\,t+c\,(c+1)\,m\,(m-1)\,t^2/2+o(t^2).
$$
On the other hand,
$$
\E((1-t)^L)=1-\E(L)\,t+\E(L(L-1))\,t^2/2+o(t^2).
$$
A comparison of the second order terms of these expansions yield the
condition on the variance of $L$ for $L$ to belong to $\ca_c$,
respectively to $\cb_c$.

To show that any square integrable distribution belongs to $\ca_c$ for
suitable values of $c$, we first 
choose values of $d$ and $s<1$, both large enough to make sure that
$\E(t^L)\le g_{d,m}(t)$ for every $t\ge s$. Thanks to the
  expansion above, this is possible for any $d$ such that 
$$
d\,(d+1)\,m\,(m-1)>\E(L(L-1)),\quad m:=\E(L).
$$ 
Then we choose a value of $c>d$
  large enough such that $1/m^{1/c}\ge g_{d,m}(s)/s$.
Thus $\E(t^L)\le g_{d,m}(t)\le g_{c,m}(t)$ for every $t\ge s$, and,
since $\E(t^L)/t$ is a nondecreasing function of $t$, for any $t\le s$,
$$
\E(t^L)\le t\,\E(s^L)/s\le t\,g_{d,m}(s)/s\le t/m^{1/c}\le g_{c,m}(t).
$$
The proof for the comparison with distributions in $\cb_c$ is similar.
\endproof
\proof{Proof of proposition~\ref{p.ex}}
Part (1) stems from the property
$$
g_{c,m}(t)\,g_{c,m'}(t)\ge t\,g_{c,mm'}(t),
$$
which we leave as an exercice for the reader.
Parts (2) and (3) are direct.
\endproof

\section{Higher moments}
\label{s.h}


Assume that $\xi_i\in\ca_c$ for every $i\le n$ and let $m:=M_n$. Then,
$$
M_n^r\,\E_k(1/Z_n^r)\le m^r\,\E_k^{c,m}(1/Z^r).
$$
Expansions of $g_{c,m}(t)$ when $m\to\infty$ show that
the distribution of $Z/m$ with respect to $\P^{c,m}_1$ converges 
to the distribution of $W$ with respect to a measure $\P_1^c$ such that
$$
\E_1^{c}(\ee^{-t\,W})=(1+c\,t)^{-1/c}.
$$
The distribution of $W$ is Gamma $(c,1/c)$, that is, its density with
respect to the Lebesgue measure $\dd w$ is
$$
w^{c-1}\,\ee^{-w/c}\,\mathbf{1}_{w\ge0}/(c^c\Gamma(c)).
$$
Furthermore, $g_{c,m}(t)\le g_c(t):=\E_1^c(t^W)=(1-c\,\log t)^{-1/c}$.
Hence,
$$
M_n^r\,\E_k(1/Z_n^r)
\le\E^c_k(1/W^r)
=
\Gamma((k/c)-r)/(c^r\,\Gamma(k/c)).
$$
This inequality holds for every positive values of $r$ and $c$ and every
positive integer $k$ such that
$k>c\,r$. The lines above prove the following result.

\begin{proposition}
\label{p.hi}
{\bf (1)}
Let $c$ such that $\xi_i\in\ca_c$ for every $i\le n$. For every
positive $r$ and $k$ such that $k>r\,c$,
$$
M_n^r\,\E_k(1/Z_n^r)\le1/[k,c]_r,
\quad
\mbox{where}\
[k,c]_r:=c^r\,\Gamma(k/c)/\Gamma((k/c)-r).
$$
When $r$ is an integer,
$$
[k,c]_r=\prod_{i=1}^r(k-ic).
$$
{\bf (2)}
Conversely, let $c$ such that $\xi_i\in\cb_c$ for every $i$.
Assume that 
$M_n\to\infty$. Then,
for every $k$ and every positive real number $r\ge k/c$,
$$
\lim_{n\to\infty}M_n^r\,\E_k(1/Z_n^r)=+\infty.
$$
\end{proposition}


\section{Bernoulli branching processes}
\label{s.g}

\subsection{Preliminaries}
We first set some notations, to be able to deal with non homogeneous processes.
\begin{definition}
The  efficiency of a Bernoulli branching 
process is the sequence 
$\cx:=(x_i)_{i\ge1}$ such that 
$\xi_i=(1-x_i)\,\delta_1+x_i\,\delta_2$.
Let $\L$, respectively $\L^*$, 
denote the set of efficiencies $\cx$ such that
$x_i\in[0,1]$, respectively $x_i\in(0,1]$, 
for every $i\ge1$.
For any $\cx\in\L$,
let $s(\cx):=(x_{i+1})_{i\ge1}$
denote the shifted sequence.\end{definition}

\begin{definition}
For any $k\ge1$ and any efficiency $\cx$, let 
$$
B_k(\cx):=\lim_{n\to\infty}\E_k(1/Z_n)\prod_{i=1}^n(1+x_i).
$$ 
In the homogeneous case
$x_i=x$ for every $i\ge1$,
we write $B_k(x)$ for
$B_k(\cx)$. 
\end{definition}

By convexity, the limit which 
defines $B_k(\cx)$ is also a supremum over $n\ge0$, thus
$B_k(\cx)\ge1/k$. 
The functional $B_k$ describes $\E_k(1/Z_n)$ for finite values of $n$
as well, since replacing every $x_i$ with $i\ge n+1$ by $0$
freezes the branching process at its value $Z_n$.
Thus, uniform upper bounds of $B_k$ on $\L$ yield upper bounds of
$\E_k(1/Z_n)$ for finite values of $n$.

\subsection{Results}
\label{s.r}

The following uniform result is a consequence of 
the fact that $\ma(\xi)<1$ for every
Bernoulli $\xi$, see proposition~\ref{p.be}.

\begin{proposition}
For every efficiency $\cx\in\L$ and every $k\ge1$, 
$$
1/k\le B_k(x)\le1/(k-1).
$$
Thus,
for every $n\ge0$,
$$
1/k\le\E_k(1/Z_n)\,\prod_{i=1}^{n}(1+x_i)\le1/(k-1).
$$
\end{proposition}

The sequence $(B_k)_{k\ge1}$ satisfies recursion relations, which
we state in proposition~\ref{p.uni}, and which characterize it fully, 
see proposition~\ref{p.exi}.

\begin{proposition}
\label{p.uni}
For every $k\ge1$, the function $B_k$ is measurable on $\L$.
Furthermore, for every $k\ge1$ and $\cx\in\L$,
\begin{equation}
\label{e.ell}
B_k(\cx)=(1+x_1)\,\sum_i{k\choose
i}\,x_1^i\,(1-x_1)^{k-i}\,B_{k+i}(s(\cx)).
\end{equation}
\end{proposition}

\begin{proposition}
\label{p.exi}
Let $(F_k)_{k\ge1}$ denote a sequence of functionals defined on $\L^*$.
Assume that
$k\,F_k(\cx)\to1$ 
when $k\to\infty$, uniformly over 
 $\cx\in\L^*$, 
and that $(F_k)_{k\ge1}$ solves (\ref{e.ell})
on $\L^*$
for every $k\ge1$. These conditions 
define a unique sequence $(F_k)_{k\ge1}$, such that
$F_k=B_k$ on $\L^*$ for every $k\ge1$.
\end{proposition}

Thus, the sequence $(B_k)_{k\ge1}$ is entirely
determined on $\L^*$ by the recursion (\ref{e.ell})
and by the bounds
$1/k\le B_k\le1/(k-1)$.

Note finally that the recursion (\ref{e.ell}) is but a special case of
the following result. For any branching process of reproducing
law $\Xi=(\xi_i)_{i\ge1}$ and any $k\ge1$, introduce
$$
H_k(\Xi):=\lim_{n\to\infty}M_n\,\E_k(1/Z_n),
$$
and
 the shifted mechanism
$s(\Xi):=(\xi_{i+1})_{i\ge1}$.
Let $\xi_1^{*k}$ denote  the convolution of the measure $\xi_1$ with
itself $k$ times.
Then,
$$
H_k(\Xi)=\mu_1\,\sum_{i\ge k}\xi_1^{*k}(i)\,H_i(s(\Xi)).
$$

\subsection{Homogeneous case}

We start with a version of the relation (\ref{e.ell}) in the
homogeneous case.

\begin{proposition}
\label{p.hh}
For every $x$ in $(0,1)$,
$$
B_k(x)=(1+x)\,\sum_i{k\choose
i}\,x^i\,(1-x)^{k-i}\,B_{k+i}(x).
$$
\end{proposition}

The recursion 
whose left hand
side is $B_k(x)$ involves the whole set of values
$B_k(x)$, $B_{k+1}(x)$,
\ldots, $B_{2k}(x)$. Thus, this system
of equations
does not yield directly
the value 
of each $B_k(x)$. The exception is the case $k=1$.

\begin{corollary}
\label{c.12}
For every
$x\neq0$,
$ 
B_1(x)=B_2(x)\,(1+x)/x.
$
\end{corollary}
Our main result in this section is proposition~\ref{p.12}.
\begin{proposition}
\label{p.12}
Let $\fb(x):=(1-x)/(1+x)$ and $\fc(x):=1-x$.
For every $k\ge1$,
$$
1/(k-\fb(x))\le
B_k(x)\le1/(k-\fc(x)).
$$
Thus, $B_k(0^+)=1/(k-1)$ 
and $B_k(1^-)=1/k$.
\end{proposition}

For $k=1$, proposition~\ref{p.12} states that $B_1(x)$ is at least
$1/(1-\alpha'(x))$.
A better bound obtains if one uses corollary~\ref{c.12} and then 
proposition~\ref{p.12}, namely
$$
(1+x)^2/(1+3x)\le x\,B_1(x)\le1.
$$
The lower bound is always greater than $8/9=.889^-$.
From our numerical simulations in section \ref{ss.sim} below, 
some values of $\lambda\, B_1(\lambda)$ are as small as
$B_*=.9274\pm.0002$, to be compared to 
$8/9=.8889^-$.

One could iterate the procedure, getting yet tighter upper and lower bounds of
$\E_1^{x}(1/Z_n)$, or of any $\E_k^{x}(1/Z_n)$
with $k\ge1$, with any prescribed accuracy.
We develop this idea in section~\ref{ss.algo} below.

We end this section with a conjecture.

\begin{problem}
We conjecture that every function $x\mapsto B_k(x)$ 
is decreasing on
$x\in(0,1]$.
Prove this and 
find a natural explanation of the fact that $B_k(0^+)$ and $B_k(0)$
are not equal.
\end{problem}

\subsection{Proofs}
\proof{Proof of proposition~\ref{p.uni}}
Since each $\E_k(1/Z_n)$ 
 is measurable with respect to $(x_i)_{i\le n}$, $B_k(x)$
 is the limit of a measurable nondecreasing sequence, hence $B_k$ is 
measurable.
As regards the recursion relation, we 
consider the conditioning by $Z_1$ of the 
Bernoulli branching process starting from $Z_0=k$.
On the event $\{Z_1=k+i\}$,
$(Z_{n+1})_{n\ge0}$ follows the law of 
a Galton-Watson branching process of efficiency
$s(x)$, starting from $k+i$. Hence (\ref{e.ell}) 
follows from the fact
that the
distribution of $Z_1-k$ is binomial $(k,x_1)$.
\endproof

\proof{Proof of proposition~\ref{p.exi}}
The existence follows from the construction of 
each $B_k$.
A proof of the uniqueness
is as follows.
Assume that the sequences of 
functionals $(B'_k)$ and $(B''_k)$ are solutions. In particular,
$B'_k/B''_k\to1$ uniformly on $\L^*$, when $k\to\infty$.
Fix $\varepsilon$.
For every $k$ large enough and for
every $x\in\L^*$,
$$
B'_k(x)\le(1+\varepsilon)\,B''_k(x).
$$
Since $(B'_k)$ and $(B''_k)$ both solve the recursion 
relations (\ref{e.ell}),
a recursion over the decreasing values of $k$ shows that
$B'_k(x)\le(1+\varepsilon)\,B''_k(x)$ for every $k\ge1$ and every
$x\in\L^*$. This recursion uses as a crucial tool the fact
that no $x_k$ is zero.
Now, since $\varepsilon$ is arbitrary, $B'_k\le B''_k$ on $\L^*$ for every
$k\ge1$.
Exchanging the roles of the two sequences, one sees that 
$B'_k=B''_k$ on $\L^*$, 
for every $k\ge1$. 
\endproof

\subsection{Outline of the proof of proposition~\ref{p.12}}
We start from
relations between the functions $B_k$ in
proposition~\ref{p.hh}, 
which read, for every $k\ge1$,
$$
B_k(x)=(1+x)\,\E_k^{x}(B_{Z_1}(x)).
$$
With respect to the probability $\P_k^{x}$, 
$Z_1$ is distributed like the sum of $k$ i.i.d.\ random variables of
distribution $(1-x)\,\delta_1+x\,\delta_2$.
Lemma \ref{l.11} below follows from the fact that $k\,B_k(x)\to1$ when
$k\to\infty$. 
\begin{lemma}
\label{l.11}
Assume that 
$\liminf k\,\varphi(k)\ge1$ and that, for every $k\ge1$,
$$
(1+x)\,\E_k^{x}(\varphi(Z_1))\le\varphi(k).
$$
Then $B_k(x)\le\varphi(k)$ for every $k\ge1$.
Conversely, if 
$\limsup k\,\psi(k)\le1$ and if, for every $k\ge1$,
$$
(1+x)\,\E_k^{x}(\psi(Z_1))\ge\psi(k),
$$
then $B_k(x)\ge\psi(k)$ for every $k\ge1$.
\end{lemma}
\begin{definition}
\label{d.22}
For every $k\ge1$, let $c_k$ denote the unique solution 
in $(0,1)$ of the
equation
$$
(1+x)\,\E_k^{x}(1/(Z_1-c_k))=1/(k-c_k).
$$
\end{definition}
Lemma \ref{l.22} becomes obvious when one uses an equivalent
definition of $c_k$, given below in part (ii) of lemma \ref{l.55}.
\begin{lemma}
\label{l.22}
If $c_k\le c$ for every $k\ge1$, then $B_k(x)\le1/(k-c)$ for
every $k\ge1$.
Conversely, if $c_k\ge c$ for every $k\ge1$, then $B_k(x)\ge1/(k-c)$ for
every $k\ge1$.
\end{lemma}
Lemma~\ref{l.22} asserts that $B_k(x)\le1/(k-c)$ for
$c=\sup\{c_k\,;\,k\ge1\}$ and, by lemma \ref{l.33}, this supremum
is $c_1=1-x$, thus proposition \ref{p.12} follows.
\begin{lemma}
\label{l.33}
For every $k\ge1$, $c_k\le c_1=1-x$.
\end{lemma}
The following result shows that the technique above cannot yield a
better value of $\fb(x)$ than
$\fb(x)=(1-x)/(1+x)$. 
\begin{lemma}
\label{l.44}
When $k\to\infty$, $c_k\to(1-x)/(1+x)$.
Furthermore, for $k\ge2$,
$$
c_k+x\,c_{k-1}\ge1-x.
$$
\end{lemma}

Finally, we use the characterizations below to evaluate 
$c_k$.

\begin{lemma}
\label{l.55}
For every $k\ge2$, the following inequalities are equivalent and
equivalent to the fact that $c\ge c_k$.
\begin{enumerate}
\item[\bf (i)]
$(1+x)\,\E_k^{x}(1/(Z_1-c))\le1/(k-c)$.
\item[\bf (ii)]
$k\,(1+x)\,\E_{k-1}^{x}(1/(Z_1+2-c))\ge1$.
\item[\bf (iii)]
$k\,(k-1)\,x\,(1+x)\,\E_{k-2}^{x}(1/(Z_1+4-c))
\le
x\,k-1+c$.
\end{enumerate}
The reversed
inequalities  {\bf (i')}, {\bf (ii')} and {\bf (iii')} are equivalent and
equivalent to the fact that $c\le c_k$.
\end{lemma}

\subsection{Technical steps of the proof of proposition~\ref{p.12}}
\label{s.p}
We prove lemmas \ref{l.33} and \ref{l.44}, assuming lemma
\ref{l.55} for the moment.
By Jensen's inequality, the expectation of the inverse is greater than
the inverse of the expectation. Thus (ii') implies
$$
k(1+x)\le(k-1)(1+x)+2-c.
$$
This reads
$c\le1-x$. Since $c_1=1-x$, we
are done with lemma \ref{l.33}.
Furthermore, we can and we will restrict the reasoning below to
$c\le1-x$.
 
To prove lemma \ref{l.44}, we first note that 
(ii) involves the
expected value of a concave
function of $u:=1/(Z_1-c_{k-1})$, namely the function
$u\mapsto u/(1+b\,u)$ with $b:=c_{k-1}+2-c$. 
The expected value of a concave function is at most the function of
the expected value. From the definition of $c_{k-1}$,
(ii)
implies
$$
k(1+x)\ge(k-1-c_{k-1})(1+x)+c_{k-1}+2-c.
$$
This is equivalent to $c\ge1-x-x\,c_{k-1}$. Hence, for
any $k\ge2$,
\begin{equation}
\label{e.cc}
1-x-x\,c_{k-1}\le c_k\le1-x.
\end{equation}
This is enough to show that $c_k\ge(1-x)^2$ for every $k\ge1$.
Thus we can and we will further restrict the 
reasoning below to $c\ge(1-x)^2$.

In the second step of the proof of lemma \ref{l.44}, we use (iii') like we
used (ii).
Namely, we note that (iii') involves the
expected value of a concave
function of $1/(Z_1-c_{k-2})$ and we apply 
Jensen's inequality once again.
From the definition of $c_{k-2}$,  (iii') implies
$$
k\,(k-1)\,x\,(1+x)\ge[(1+x)\,(k-2-c_{k-2})+c_{k-2}+4-c]\,
(x\,k-1+c).
$$
After some simplifications, this reads $A_1\,k+A_0\ge0$, with
$$
A_1:=(1-x)^2-c+x^2\,c_{k-2},
\quad
A_0:=(1-c)(2(1-x)-c-x\,c_{k-2}).
$$
Since $c\ge(1-x)^2$ and $c_{k-2}\ge(1-x)^2$, simple bounds
show that $A_0\le1$. 
Hence (iii') implies that $A_1\ge-1/k$. Finally, for every $k\ge3$,
\begin{equation}
\label{e.ccc}
(1-x)^2\le c_k\le(1-x)^2+x^2\,c_{k-2}+1/k.
\end{equation}
One uses the a priori bounds of (\ref{e.cc}) and (\ref{e.ccc}) as follows.
On the one hand, the upper bound of $c_k$ in (\ref{e.ccc}) implies
$$
\limsup c_k\le(1-x)^2+x^2\,\limsup c_{k}.
$$
On the other hand, the lower bound of $c_k$ in
(\ref{e.cc}) implies
$$
1-x-x\,\limsup c_k\le
\liminf c_k.
$$
Hence $\limsup c_k=\liminf c_k=(1-x)/(1+x)$.
This proves lemma \ref{l.44}.

Lemma~\ref{l.55} is a consequence of the following trick.
Part (i) involves 
$$
(k-c)/(Z_1-c)=1-(Z_1-k)/(Z_1-c)=:1-v.
$$
By exchangeability, $\E_k^{x}(v)$  is $k$
times the expected value of $(L_1-1)/(Z_1-c)$, where $L_1$ denotes the
number of descendants of the first individual in the initial
population. The event $\{L_1-1\neq0\}$ is $\{L_1=2\}$ and 
has probability $x$. Thus, for every $k\ge2$,
$$
(k-c)\,\E^{x}_k(1/(Z_1-c))=1-k\,x\,\E^{x}_{k-1}(1/(Z_1+2-c)).
$$
With the convention that $Z_1=0$, $\P^{x}_0$ almost surely, this
relation holds for $k=1$ as well.
This translates (i) or (i') into (ii) or (ii').
The translation of (ii) or (ii') into (iii) or (iii') uses the same
trick, starting from $1/(Z_1+2-c)$. This concludes the proof of
proposition~\ref{p.12}.

\subsection{Algorithm}
\label{ss.algo}

The  following algorithm
yields approximate values of $B_k$ on $\L^*$, 
with any prescribed accuracy. 
\begin{itemize}
\item
Fix $n\ge1$. 
\item
For every $k\ge n+1$ and $x$, 
let 
$$
B^0_{k,n}(x):=1/k,
\quad
B^1_{k,n}(x):=1/(k-1). 
$$
\item
Find the unique sequence
$(B^1_{k,n})_{k\le n}$ that solves the system of equations
(\ref{e.ell}) for $k\le n$, when one replaces every 
$B_k(s(x))$ such that
$k\ge n+1$, by $B^1_{k,n}(s(x))$, that is, by the value $1/(k-1)$.
\item
Likewise, find the unique sequence
$(B^0_{k,n}(x))_{k\le n}$ that solves the system of equations
(\ref{e.ell}) for $k\le n$, when one replaces every 
$B_k(s(x))$ such that
$k\ge n+1$, by $B^1_{k,n}(x)$, that is, by the value $1/k$.
\item
Then, for every $k\ge1$ and every $x$,
$$
B^0_{k,n}(x)\le B_k(x)\le B^1_{k,n}(x)
\le (1+1/n)\,B^0_{k,n}(x).
$$
\end{itemize}

\subsection{Comments on the algorithm}

Neither $(B^0_{k,n})_{k\ge1}$ nor $(B^1_{k,n})_{k\ge1}$
solve the full system of equations
(\ref{e.ell}).
For any fixed values of $k$ and $x$, 
$(B^0_{k,n}(x))_{n\ge1}$ is a nondecreasing
sequence that converges to $B_k(x)$ when $n\to\infty$,
Likewise, $(B^1_{k,n}(x))_{n\ge1}$ is a nonincreasing
sequence that converges to $B_k(x)$ when $n\to\infty$,

Increasing values of $n$ yield more and more accurate approximations
of each $B_k(x)$ and the relative error 
is of order at most $1/n$.

In the Bernoulli case, one can use some initial values, better than
$1/k$, respectively $1/(k-1)$, namely, for every $k\ge n+1$ and $x$,  
$$
b^0_{k,n}(x):=1/(k-\fb(x)),
\quad
b^1_{k,n}(x):=1/(k-\fc(x)). 
$$
The relative error that was at most $1+1/n$ in the first version of
the algorithm becomes
at most 
$$
1+(\fc-\fb)(x)/(n+1-\fc(x))\le1+(3-2\sqrt{2})/(n+x).
$$
Numerically, this is at most $1+.172/n$, for every $x$.

\subsection{Simulations in the homogeneous case}
\label{ss.sim}
The algorithm above 
with $n:=1000$ suggest the following refinements.
Define 
$$
B(x):=B_1(x)\,x=B_2(x)\,(1+x).
$$
Simulations show that
$B$ decreases on $(0,x_*)$
from $B(0^+)=1$ to $B(x_*):=B_*$
and increases on $(x_*,1]$ from $B_*$ to
  $B(1)=1$, with
$$
x_*=.38\pm.01,
\quad
B_*=.9274\pm.0002.
$$
This would imply that, 
for every positive $x$, 
$$
B_*/x\le B_1(x)\le1/x.
$$
Simulations show that $B_2$, hence $B_1$ as well, 
decreases on $(0,1]$.
%


\section{Size-biased offspring}
\label{s.sbo}
When computing harmonic means, it may prove convenient to use size-biased
distributions, as follows. 
Assume that $L$ and $L_i$ are i.i.d.\ positive integrable
random variables and
that $L'$ is an independent size-biased copy of $L$, that is, for
every $t\in[0,1]$,
$$
\E(t^{L'}):=\E(L\,t^L)/\E(L).
$$
Then,  for any nonnegative integer $k$,
$$
\E(L)\,\E(1/(L_1+\cdots+L_k+L'))=1/(k+1).
$$
Can one use this in our branching setting?
Assume first that $1\le L\le c+1$ almost surely, for a given integer $c$. 
Since $L'\le c+1\le L_{k+1}+\cdots+L_{k+c+1}$ almost surely, this
proves that
$$
\E(L)\,\E_{k}(1/Z_1)\le1/(k-c),
$$
for every $k\ge c+1$.
More generally, the last inequality above
 holds as soon as
$\E(t^{L'})\ge\E(t^L)^{c+1}$ for every $t\in[0,1]$ and for a given
positive $c$.

In our setting, this line of reasoning suffers from two
drawbacks. First, unless we miss something, to be able to
iterate this inequality during $n$
generations, one must assume that $k\ge
n\,c$.
Second, the inequality $\E(t^{L'})\ge\E(t^L)^{c+1}$ implies that
$c\ge\ma(\xi)$, where $\xi$ denotes the law of $L$ (the proof is easy
and omitted). In other words, 
$k\ge c$
implies that $\xi\in\ca_c$.

\section{Case $k\le c$}
\label{s.cri}

This section is a brief description of 
the behaviour of $\E_k(1/Z_n)$ when the
hypotheses of theorem~\ref{t.b} are not met.
Consider, for the sake of simplicity, a homogeneous
branching process and let 
$$
p_i:=\xi(i)=\P(L=i)=\P(Z_1=i\,\vert\,Z_0=1).
$$
Our first remark is that, for every $n\ge0$ and $k\ge1$,
$$
\E_k(1/Z_n)\ge r_k^n/k,\quad
r_k:=\max\{p_1^k,1/\mu\}.
$$
The $1/\mu$ bound is 
due to the convexity.
The $p_1^k$ bound is due to the fact that
the probability of the event $\{Z_n=k\}$ is $p_1^k$.
 
The parameters $\mu$ and $p_1$, through $r_k$, indeed 
describe the asymptotics of
$\E_k(1/Z_n)$, as follows.
For the sake of simplicity, we exclude
the degenerate case $p_1^k\,\mu=1$,
where polynomial corrections appear.
For a given $k\ge1$, 
there exists a finite positive $h_k$ such that
$$
\lim_{n\to\infty}\E_k(1/Z_n)/r_k^n=h_k.
$$
In the Bernoulli case, $p_1<1/\mu$ hence
$r_k=1/\mu$ for every $k\ge1$ and the $r_k=p_1^k$
regime is nonexistent.

The limits $h_k$ satisfy the following relations. Assume for instance
that one wishes to compute $h_1$ and that
$\mu\,p_1>1$, hence $r_1=p_1$.
Conditioning on
the value of
$Z_1$, one gets a relation between $\E_1(1/Z_{n+1})$ and 
the sequence $(\E_k(1/Z_n))_{k\ge1}$.
Letting $n$ go to infinity yields
$$
h_1=1+\sum_{k\ge2}H_k\,p_k/p_1, \quad
H_k:=\sum_{n\ge0}\E_k(1/Z_n)/p_1^{n}.
$$
The term $\E_2(1/Z_n)/p_1^n$
behaves like $(r_2/p_1)^n$, that is, like
$1/(\mu\,p_1)^n$ if
$\mu\,p_1^2>1$,
or like $p_1^n$ if $\mu\,p_1^2<1$.
Since both quantities are summable,
$H_2$ is finite. Since $H_k\le H_2$ for every $k\ge2$,
$h_1$ is finite as well.

When $\mu\,p_1<1$, $r_k=1/\mu$ for every $k\ge1$ and
the same reasoning as above yields
$$
h_1=\mu\,\sum_{k\ge1}p_k\,h_k.
$$
Since $\mu\,p_1<1$, this gives $h_1$ as a linear combination of
$(h_k)_{k\ge2}$.
In turn, for every  $k\ge2$, a one-step recursion, similar to 
the one we used before, shows that $h_k$ is such that
$$
h_k=\mu\,\sum_{i\ge k}h_i\,\xi^{*k}(i).
$$
Since the
coefficient of $h_k$ in the right hand side is $\mu\,p_1^k<1$, 
this relation implies that $h_k$ is a linear combination of
the sequence $(h_i)_{i\ge k+1}$ and that this linear combination uses
nonnegative coefficients. 
Since $h_i\le h_2$ for every $i\ge k+1$,
this series converges.
However, it does not seem easy to get information about the
coefficients $h_1$ or $h_k$
with $k\ge2$ from these relations.




\bigskip\bigskip

\begin{tabular}{l}
Universit\'e Claude Bernard Lyon 1 \\ 
Institut Camille Jordan UMR 5208 
\\
Domaine de Gerland \\
50, avenue Tony-Garnier \\ 69366 Lyon Cedex 07 (France) \\ \ \\ {\tt
Didier.Piau@univ-lyon1.fr} \\ {\tt http://lapcs.univ-lyon1.fr}
\end{tabular}

\end{document}